\documentclass[12pt,eqno,]{article}
\usepackage{amsmath,amsthm,amssymb}

\makeatletter

\newcommand{\Rmnum}[1]{\expandafter\@slowromancap\romannumeral #1@}
\makeatother

\newcommand{\arth}{\mbox{arctanh}}

\title{Sharp two parameter bounds for logarithmic and arithmetic-geometric means}
\author{\small Yu-Ming Chu$^{1}$, Ye-Fang Qiu$^{1}$, Miao-Kun Wang$^{1}$ and Xiao-Yan Ma$^{2}$}
\date{}

\begin{document}
\maketitle

{\footnotesize\rm
\noindent $^{1}$Department of Mathematics, Huzhou Teachers College, Huzhou 313000, China;\\
$^{2}$Department of Mathematics, Zhejiang Sci-Tech University, Hangzhou 310018, China.\\
Correspondence should be addressed to Yu-Ming Chu,
chuyuming@hutc.zj.cn}

\bigskip \noindent{\bf Abstract}: For fixed $s\geq 1$ and
$t_{1},t_{2}\in(0,1/2)$ we prove that the inequalities
$G^{s}(t_{1}a+(1-t_{1})b,t_{1}b+(1-t_{1})a)A^{1-s}(a,b)>AG(a,b)$ and
$G^{s}(t_{2}a+(1-t_{2})b,t_{2}b+(1-t_{2})a)A^{1-s}(a,b)>L(a,b)$ hold
for all $a,b>0$ with $a\neq b$ if and only if $t_{1}\geq
1/2-\sqrt{2s}/(4s)$ and $t_{2}\geq 1/2-\sqrt{6s}/(6s)$. Here
$G(a,b)$, $L(a,b)$, $AG(a,b)$ and $A(a,b)$ are the geometric,
logarithmic, arithmetic-geometric and arithmetic means of $a$ and
$b$, respectively.



\section{Introduction}
\hspace{4mm} For real numbers $a,b$ and $c$ with $c\neq
0,-1,-2,\cdots$, the Gaussian hypergeometric function is defined by
\begin{equation}
F(a,b;c;x)={}_{2}F_{1}(a,b;c;x)=\sum_
{n=0}^{\infty}\frac{(a,n)(b,n)}{(c,n)}\frac{x^n}{n!},\ |x|<1.
\end{equation}
Here $(a,0)=1$ for $a\neq0$, and
$(a,n)=a(a+1)(a+2)(a+3)\cdots(a+n-1)$ is the shifted factorial
function for $n=1,2,\cdots$. In connection with the Gaussian
hypergeometric function, the well-known complete elliptic integrals
$\mathcal{K}(r)$ and $\mathcal{E}(r)(0<r<1)$ of the first and second
kinds [1, 2] are defined by
\begin{equation}
\left\{\begin{array}{ll}
{\mathcal{K}}(r)={\pi}F({1}/{2},{1}/{2};1;r^2)/2=\int_{0}^{{\pi}/{2}}(1-r^2\sin^2\theta)^{-1/2}d \theta,\\
{\mathcal{K}}(0)={\pi}/2,\qquad{\mathcal{K}}(1)=\infty
\end{array}\right.
\end{equation}
and
\begin{equation}
\left\{\begin{array}{ll}
{\mathcal{E}}(r)={\pi}F(-{1}/{2},{1}/{2};1;r^2)/2=\int_{0}^{{\pi}/{2}}(1-r^2\sin^2\theta)^{1/2} d\theta,\\
{\mathcal{E}}(0)={\pi}/2,\qquad{\mathcal{E}}(1)=1,
\end{array}\right.
\end{equation}
respectively. The following formulas for $\mathcal{K}(r)$ were
presented in [3]:
\begin{equation}
\frac{d{\mathcal{K}(r)}}{dr}=\frac{{\mathcal{E}(r)}-{(1-r^2)}{\mathcal{K}(r)}}{r(1-r^2)},
\end{equation}
\begin{equation}
\mathcal{K}\left(\frac{2\sqrt{r}}{1+r}\right)=(1+r)\mathcal{K}(r).
\end{equation}

Let $H(a,b)=2ab/(a+b)$, $G(a,b)=\sqrt{ab}$,
$L(a,b)=(b-a)/(\log{a}-\log{b})$ and $A(a,b)=(a+b)/2$ be the
classical harmonic, geometric, logarithmic and arithmetic means of
two distinct positive real numbers $a$ and $b$, respectively. Then
it is well known that the inequalities $H(a,b)<G(a,b)<L(a,b)<A(a,b)$
hold for all $a,b>0$ with $a\neq b$.

The classical arithmetic-geometric mean $AG(a,b)$ of two positive
number $a$ and $b$ is defined as the common limit of sequences
$\{a_{n}\}$ and $\{b_{n}\}$, which are given by
\begin{align*}
&a_{0}=a,\qquad b_{0}=b,\\
&a_{n+1}=(a_{n}+b_{n})/2=A(a_{n},b_{n}),\qquad
b_{n+1}=\sqrt{a_{n}b_{n}}=G(a_{n},b_{n}).
\end{align*}

It is well known that inequalities
\begin{equation}
G(a,b)<\sqrt{A(a,b)G(a,b)}<AG(a,b)<A(a,b)
\end{equation}
hold for all $a,b>0$ with $a\neq b$.

Recently, the harmonic, geometric, logarithmic, arithmetic-geometric
and arithmetic means have been the subject of intensive research. In
particular, many remarkable inequalities for these means can be
found in the literature [4-13].

The Gaussian identity [3] shows that
\begin{equation}
AG(1,r)\mathcal{K}(\sqrt{1-r^2})=\frac{\pi}{2}
\end{equation}
for all $r\in(0,1)$.

Carlson and Vuorinen [14], and Brackenn [15] proved that
\begin{equation*}
L(a,b)<AG(a,b)
\end{equation*}
for all $a,b>0$ with $a\neq b$. Vamanamurthy and Vuorinen [16]
established that
\begin{equation*}
AG(a,b)<\frac{\pi}{2}L(a,b)
\end{equation*}
for all $a,b>0$ with $a\neq b$.

For $t_{1},t_{2},t_{3},t_{4}\in(0,1/2)$, very recently Chu et al.
[17, 18] proved that the inequalities
\begin{equation}
G(t_{1}a+(1-t_{1})b,t_{1}b+(1-t_{1})a)>AG(a,b),
\end{equation}
\begin{equation}
H(t_{2}a+(1-t_{2})b,t_{2}b+(1-t_{2})a)>AG(a,b),
\end{equation}
\begin{equation}
G(t_{3}a+(1-t_{3})b,t_{3}b+(1-t_{3})a)>L(a,b),
\end{equation}
and
\begin{equation}
H(t_{4}a+(1-t_{4})b,t_{4}b+(1-t_{4})a)>L(a,b)
\end{equation}
hold for all $a,b>0$ with $a\neq b$ if and only if $t_{1}\geq
1/2-\sqrt{2}/4$, $t_{2}\geq 1/4$, $t_{3}\geq 1/2-\sqrt{6}/6$ and
$t_{4}\geq 1/2-\sqrt{3}/6$.

Let $t\in(0,1/2)$, $s\geq 1$ and
\begin{equation}
Q_{t,s}(a,b)=G^{s}(ta+(1-t)b, tb+(1-t)a)A^{1-s}(a,b).
\end{equation}
Then it is not difficult to verify that
\begin{equation*}
Q_{t,1}(a,b)=G(ta+(1-t)b, tb+(1-t)a),
\end{equation*}
\begin{equation*}
Q_{t,2}(a,b)=H(ta+(1-t)b, tb+(1-t)a)
\end{equation*}
and $Q_{t,s}(a,b)$ is strictly increasing with respect to
$t\in(0,1/2)$ for fixed $a,b>0$ with $a\neq b$.

It is natural to ask what are the least values $t_{1}=t_{1}(s)$ and
$t_{2}=t_{2}(s)$ in $(0,1/2)$ such that inequalities
$Q_{t_{1},s}(a,b)>AG(a,b)$ and $Q_{t_{2},s}(a,b)>L(a,b)$ hold for
all $a,b>0$ with $a\neq b$ and $s\geq 1$. The aim of this paper is
to answer these questions, our main results are the following
Theorems 1.1 and 1.2.

\medskip
{\bf Theorem 1.1.} If $t\in(0,1/2)$ and $s\geq 1$ then the
inequality
\begin{equation}
Q_{t,s}(a,b)>AG(a,b)
\end{equation}
holds for all $a,b>0$ with $a\neq b$ if and only if $t\geq
1/2-\sqrt{2s}/(4s)$.

\medskip
{\bf Theorem 1.2.} If $t\in(0,1/2)$ and $s\geq 1$ then the
inequality
\begin{equation}
Q_{t,s}(a,b)>L(a,b)
\end{equation}
holds for all $a,b>0$ with $a\neq b$ if and only if $t\geq
1/2-\sqrt{6s}/(6s)$.

\medskip
{\bf Remark 1.1.} Let $s=1,2$ in Theorem 1.1, then inequality (1.13)
becomes inequalities (1.8) and (1.9), respectively.

\medskip
{\bf Remark 1.2.} Let $s=1,2$ in Theorem 1.2, then inequality (1.14)
becomes inequalities (1.10) and (1.11), respectively.

\section{Lemmas}
\hspace{4mm} \setcounter{equation}{0}

In order to prove Theorems 1.1 and 1.2 we need two lemmas, which we
present in this section.

\medskip
{\bf Lemma 2.1.} Let $u\in[0,1]$, $s\geq 1$ and
\begin{equation}
f_{u,s}(x)=\frac{s}{2}\log(1-ux^2)-\log\left(\frac{\pi}{2\mathcal{K}(x)}\right).
\end{equation}
Then $f_{u,s}>0$ for all $x\in(0,1)$ if and only if $2su\leq 1$.

\medskip
{\bf\em Proof.} From (1.4) and (2.1) one has
\begin{equation}
f_{u,s}'(x)=-\frac{usx}{1-ux^2}+\frac{\mathcal{E}(x)-(1-x^2)\mathcal{K}(x)}{x(1-x^2)\mathcal{K}(x)}=\frac{F_{u,s}(x)}{x(1-x^2)(1-ux^2)\mathcal{K}(x)},
\end{equation}
where
\begin{equation}
F_{u,s}(x)=-sux^2(1-x^2)\mathcal{K}(x)+(1-ux^2)[\mathcal{E}(x)-(1-x^2)\mathcal{K}(x)].
\end{equation}

It follows from (1.1)-(1.3) and (2.3) together with elaborated
computations that
\begin{align*}
&\mathcal{E}(x)-(1-x^2)\mathcal{K}(x)\\
=&\frac{\pi}{2}\left[\sum\limits_{n=0}^{\infty}\frac{(-1/2,n)(1/2,n)}{(n!)^2}x^{2n}-(1-x^2)\sum\limits_{n=0}^{\infty}\frac{(1/2,n)^2}{(n!)^2}x^{2n}\right]\\
=&\frac{\pi}{2}\sum\limits_{n=0}^{\infty}\frac{(1/2,n)^2}{2n!(n+1)!}x^{2n+2},
\end{align*}
\begin{align}
\frac{2}{\pi}F_{u,s}(x)=&-sux^2(1-x^2)\sum\limits_{n=0}^{\infty}\frac{(1/2,n)^2}{(n!)^2}x^{2n}+(1-ux^2)\sum\limits_{n=0}^{\infty}\frac{(1/2,n)^2}{2n!(n+1)!}x^{2n+2}\nonumber\\
=&-su\sum\limits_{n=0}^{\infty}\frac{(1/2,n)^2}{(n!)^2}x^{2n+2}+su\sum\limits_{n=0}^{\infty}\frac{(1/2,n)^2}{(n!)^2}x^{2n+4}\nonumber\\
&+\sum\limits_{n=0}^{\infty}\frac{(1/2,n)^2}{2n!(n+1)!}x^{2n+2}-u\sum\limits_{n=0}^{\infty}\frac{(1/2,n)^2}{2n!(n+1)!}x^{2n+4}\nonumber\\
=&-sux^2-su\sum\limits_{n=0}^{\infty}\frac{(1/2,n+1)^2}{[(n+1)!]^2}x^{2n+4}+su\sum\limits_{n=0}^{\infty}\frac{(1/2,n)^2}{(n!)^2}x^{2n+4}\nonumber\\
&+\frac{x^2}{2}+\sum\limits_{n=0}^{\infty}\frac{(1/2,n+1)^2}{2(n+1)!(n+2)!}x^{2n+4}-u\sum\limits_{n=0}^{\infty}\frac{(1/2,n)^2}{2n!(n+1)!}x^{2n+4}\nonumber\\
=&x^2\left[\frac{1}{2}-su+\sum\limits_{n=0}^{\infty}\frac{(1/2,n)^2A_{n}}{2(n+1)!(n+2)!}x^{2n+2}\right],
\end{align}
where
\begin{equation}
A_{n}=su(n+2)(2n+\frac{3}{2})+(n+\frac{1}{2})^2-u(n+1)(n+2)>0.
\end{equation}

We divide the proof into two cases.

{\bf\em Case 1.1.} $2su\leq 1$. Then (2.2)-(2.5) lead to conclusion
that $f_{u,s}(x)$ is strictly increasing in $(0,1)$. Therefore,
$f_{u,s}(x)>f_{u,s}(0^+)=0$ for all $x\in(0,1)$ follows from (1.2)
and (2.1) together with the monotonicity of $f_{u,s}(x)$ in $(0,1)$.

{\bf\em Case 1.2.} $2su>1$. Then (2.2)-(2.4) lead to conclusion that
there exists $\delta_{1}\in(0,1)$ such that $f_{u,s}(x)$ is strictly
decreasing in $(0,\delta_{1})$. Therefore,
$f_{u,s}(x)<f_{u,s}(0^+)=0$ for all $x\in(0,\delta_{1})$ follows
from (1.2) and (2.1) together with the monotonicity of $f_{u,s}(x)$
in $(0,\delta_{1})$.

\medskip
{\bf Lemma 2.2.} Let $u\in [0,1]$, $s\geq 1$,
$\arth(x)=\frac{1}{2}\log\left(\frac{1+x}{1-x}\right)$ be the
inverse hyperbolic tangent function, and
\begin{equation}
g_{u,s}(x)=\frac{s}{2}\log(1-ux^2)+\log\left(\frac{\arth(x)}{x}\right).
\end{equation}
Then $g_{u,s}(x)>0$ for all $x\in(0,1)$ if and only if $3su\leq 2$.

\medskip
{\bf\em Proof.} From (2.6) one has
\begin{equation}
g_{u,s}'(x)=-\frac{sux}{1-ux^2}+\frac{x-(1-x^2)\arth{(x)}}{x(1-x^2)\arth{(x)}}=\frac{G_{u,s}(x)}{x(1-x^2)(1-ux^2)\arth(x)},
\end{equation}
where
\begin{equation}
G_{u,s}(x)=-sux^2(1-x^2)\arth(x)+(1-ux^2)[x-(1-x^2)\arth(x)].
\end{equation}

Making use of series expansion and (2.8) we have
\begin{align}
G_{u,s}(x)=&-sux^2(1-x^2)\sum\limits_{n=0}^{\infty}\frac{x^{2n+1}}{2n+1}+(1-ux^2)\left[x-(1-x^2)\sum\limits_{n=0}^{\infty}\frac{x^{2n+1}}{2n+1}\right]\nonumber\\
=&-su\sum\limits_{n=0}^{\infty}\frac{x^{2n+3}}{2n+1}+su\sum\limits_{n=0}^{\infty}\frac{x^{2n+5}}{2n+1}+(1-ux^2)
\sum\limits_{n=0}^{\infty}\frac{2x^{2n+3}}{(2n+1)(2n+3)}\nonumber\\
=&x^3\left[\frac{2}{3}-su+\sum\limits_{n=0}^{\infty}\frac{B_{n}x^{2n+2}}{(2n+1)(2n+3)(2n+5)}\right],
\end{align}
where
\begin{equation}
B_{n}=2u(s-1)(2n+5)+2(2n+1)>0.
\end{equation}

We divide the proof into two cases.

{\bf\em Case 1.1.} $3su\leq 2$. Then (2.7)-(2.10) lead to conclusion
that $g_{u,s}(x)$ is strictly increasing in $(0,1)$. Therefore,
$g_{u,s}(x)>g_{u,s}(0^+)=0$ for all $x\in(0,1)$ follows from (2.6)
together with the monotonicity of $g_{u,s}(x)$ in $(0,1)$.

{\bf\em Case 1.2.} $3su>2$. Then (2.7)-(2.9) lead to conclusion that
there exists $\delta_{2}\in(0,1)$ such that $g_{u,s}(x)$ is strictly
decreasing in $(0,\delta_{2})$. Therefore,
$g_{u,s}(x)<g_{u,s}(0^+)=0$ for all $x\in(0,\delta_{2})$ follows
from (2.6) and the monotonicity of $g_{u,s}(x)$ in $(0,\delta_{2})$.

\section{Proof of Theorems 1.1 and 1.2}
\hspace{4mm} \setcounter{equation}{0} {\bf Proof of Theorem 1.1.}
Since both $Q_{t,s}(a,b)$ and $AG(a,b)$ are symmetric and
homogeneous of degree 1. Without loss of generality, we assume that
$a>b$. Let $x=(a-b)/(a+b)\in(0,1)$. Then from (1.5) and (1.7)
together with $b/a=(1-x)/(1+x)$ we have
\begin{align}
\frac{AG(a,b)}{A(a,b)}=&\frac{AG(1,b/a)}{A(1,b/a)}=\frac{\pi}{\mathcal{K}\sqrt{1-(b/a)^2}(1+b/a)}\nonumber\\
=&\frac{\pi(1+x)}{2\mathcal{K}(2\sqrt{x}/(1+x))}=\frac{\pi}{2\mathcal{K}(x)}.
\end{align}

It follow from (1.12) and (3.1) that
\begin{align}
&\log\left(\frac{Q_{t,s}(a,b)}{AG(a,b)}\right)=\log\left(\frac{Q_{t,s}(a,b)}{A(a,b)}\right)-\log\left(\frac{AG(a,b)}{A(a,b)}\right)\nonumber\\
&=\frac{s}{2}\log\left[1-(1-2t)^2x^2\right]-\log\left[\frac{\pi}{2\mathcal{K}(x)}\right].
\end{align}

Therefore, Theorem 1.1  follows from Lemma 2.1 and (3.2).

\medskip
{\bf Proof of Theorem 1.2.} Since both $Q_{t,s}(a,b)$ and $L(a,b)$
are symmetric and homogeneous of degree 1. Without loss of
generality, we assume that $a>b$. Let $x=(a-b)/(a+b)\in(0,1)$. Then
(1.12) leads to
\begin{align}
&\log\left(\frac{Q_{t,s}(a,b)}{L(a,b)}\right)=\log\left(\frac{Q_{t,s}(a,b)}{A(a,b)}\right)-\log\left(\frac{L(a,b)}{A(a,b)}\right)\nonumber\\
&=\frac{s}{2}\log\left[1-(1-2t)^2x^2\right]+\log\left(\frac{\arth(x)}{x}\right).
\end{align}

Therefore, Theorem 1.2 follows from Lemma 2.2 and (3.3).\\

{\bf Acknowledgement:} This work was supported by the Natural
Science Foundation of China (Grant Nos. 11071059, 11071069,
11171307), and the Innovation Team Foundation of the Department of
Education of Zhejiang Province (Grant no. T200924).

\end{document}